\documentclass[10pt]{article}

\usepackage{amssymb}
\usepackage{amsmath}
\usepackage{float}
\usepackage{verbatim}
\usepackage{url}
\usepackage{graphicx}
\usepackage{subfigure}
\usepackage[dvips]{epsfig}

\newcommand{\R}{\mbox{${\mathbb R}$}}
\newcommand{\cD} {\mbox{$\cal D$}}

\newcommand{\ip}[2]{\left\langle #1,#2 \right\rangle}

\newtheorem{theorem}{Theorem}[section]

\newtheorem{algorithm}[theorem]{Algorithm}

\begin{document}

\author{Steven J. Benson\footnote{Mathematics and Computer Science Division, 
Argonne National Laboratory, Argonne, Illinois 60439, 
\texttt{benson@mcs.anl.gov}}
\and Todd S. Munson\footnote{
Mathematics and Computer Science Division, 
Argonne National Laboratory, Argonne, Illinois 60439, 
\texttt{tmunson@mcs.anl.gov}}} 

\title{Flexible Complementarity Solvers for\\Large-Scale Applications\thanks{This work was supported by the Mathematical, Information, and Computational Sciences Division subprogram of the Office of Advanced Scientific Computing, Office of Science, U.S. Department of Energy, under Contract W-31-109-Eng-38.}}  

\date{Submitted: \today} 

\maketitle


\begin{abstract}
  Discretizations of infinite-dimensional variational inequalities
  lead to linear and nonlinear complementarity problems with many
  degrees of freedom.  To solve these problems in a parallel computing
  environment, we propose two active-set methods that solve only one linear
  system of equations per iteration.  The linear solver,
  preconditioner, and matrix structures can be chosen by the user for
  a particular application to achieve high parallel performance.
  The parallel scalability of these methods is demonstrated for some
  discretizations of infinite-dimensional variational inequalities.
\end{abstract}

\section{Introduction}


Achieving high performance for numerical methods in parallel computing
environments demands that the user have the ability to customize
algorithms, linear solvers, and data structures for their particular
problems.  Closed environments or algorithms requiring specific linear
algebra constrain the choices available to the user, inevitably
leading to inefficiency.  This paper concerns a flexible, open
environment, the Toolkit for Advanced Optimization, and two algorithms
for solving complementarity problems implemented so that the user can
exploit problem structure.  The benefits of this design are
significant reductions in solution time and good
parallel scalability on targeted complementarity problems.


Complementarity problems arise in various application areas
(see, e.g., \cite{harker.pang:finite-dimensional,ferris.pang:engineering}).
In this paper, we are interested particularly in applications arising
from infinite-dimensional box-constrained variational inequalities 
where a discretization of the problem corresponds to a large-scale
linear or nonlinear complementarity problem
\cite{ekeland.temam:convex,kinderlehrer.stampacchia:introduction}.  One 
example of this
type, the journal bearing problem \cite{averick.carter.ea:minpack-2},
involves determining the pressure distribution of a thin film of lubricant
between two circular cylinders.  This problem is posed as an elliptic
partial differential equation with a free boundary.  A
finite-difference scheme for solving instances of this problem has
$n_xn_y$ degrees of freedom, where $n_x$ and $n_y$ are the number of
points in the spatial discretizations.  The large number of degrees of
freedom and structure imposed by the finite-difference scheme implies
that high performance can be achieved when solving this problem in a
parallel environment.  Other applications that can be posed as
infinite-dimensional variational inequalities include pricing American
options \cite{huang.pang:option,wilmott.dewynne.ea:option}, nonlinear
obstacle problems \cite{rodrigues:obstacle}, and optimal control
problems \cite{ulbrich:nonsmooth}.



The Toolkit for Advanced Optimization (TAO) provides a flexible
environment for solving optimization and complementarity problems in
parallel \cite{tao-user-ref}.  Several numerical linear algebra
packages have been incorporated into TAO
\cite{balay.gropp.ea:petsc,hypre,trilinos-web-page}.  These packages
offer many algorithmic choices so that a user can select the most
appropriate method for a particular application.  Novice users can
rely on the provided defaults, while more advanced users have the
freedom to select their problem representation and linear algebra to
improve parallel performance.  All of the source code for TAO can be
downloaded from \cite{tao-web-page}, which includes our
implementations of several algorithms for solving complementarity
problems.  


Section 2 discusses the design philosophy of and facilities
provided by the Toolkit for Advanced Optimization.
The complementarity algorithms used to solve the discretized nonlinear
complementarity problems in a parallel computing environment are
discussed in Section 3.  Two methods were implemented and tested: a
semismooth method
\cite{deluca.facchinei.ea:semismooth,munson.facchinei.ea:semismooth}
and a reduced space method similar to those presented in
\cite{dirkse.ferris:crash,facchinei.soares:merit}.  These algorithms
are attractive because they solve only a single system of linear
equations per iteration.  These solves can be performed by using
preconditioned iterative methods \cite{saad:iterative}, taking
advantage of the research lavished on numerical linear algebra for
solving partial differential equations.  Numerical results on the
MCPLIB collection \cite{dirkse.ferris:mcplib} of complementarity
problems are presented in Section 4 to demonstrate that the methods
are reasonably robust on general problems.  The parallel
performance achieved on several discretizations of
infinite-dimensional variational inequalities shows the
benefits of customizing the linear algebra:
both good parallel performance and scalability are achieved.


\section{Design Philosophy}

Traditionally, the numerical methods used in an algorithm implementation
to solve linear systems of equations are chosen by the developer.  
Since most floating point
operations are typically performed by the linear system solver,
considerable effort has been expended to select and implement one that
is efficient and robust for a diverse collection of applications.
Efficiency for particular problem instances has been sacrificed in exchange
for robustness and general applicability.  However, linear solvers
tailored to a particular application may save a significant amount of
computation and allow the method to find solutions to larger problems
than previously possible.  While developers choose a linear solver
without any knowledge of the application, many applications have
structure that can be exploited.  The Toolkit for Advanced
Optimization (TAO) \cite{tao-web-page,Benson:2001:CSP,tao-user-ref}
was specifically designed to allow the user to customize the linear
solver to their application so that they can achieve high performance
and parallel scalability.

One linear solver choice is an LU factorization, which can be
applied to arbitrary nonsingular systems.  Many robust
implementations have been developed that utilize sparsity 
while maintaining numerical
stability.  The complementarity methods in TAO
allow for the use of the LU factorizations from LAPACK
\cite{anderson.bai.ea:lapack} for dense systems and LUSOL
\cite{murtagh.saunders:minos} for sparse systems.  Nonetheless, these
factorizations do not exploit symmetry, strong monotonicity, or other
features of the matrix, such as block structure.  Furthermore, direct
factorizations can impose excessive memory requirements even when the
given matrix is sparse.  Access to Cholesky factorizations are provided
when a symmetric positive definite linear system is solved.  In this
case, the number of floating-point operations can be reduced by a
factor of two, and further gains can be achieved by stable reorderings
of the matrix based solely on the sparsity pattern.

Aside from direct methods, many iterative techniques can be applied to
solve linear systems of equations.  These techniques typically do not
impose the memory requirements of direct methods.  Two of the most
common techniques are GMRES \cite{saad.schultz:gmres} for general
nonsingular matrices and the conjugate gradient method
\cite{golub.vanloan:matrix} for symmetric positive definite matrices.
Convergence of these methods can be significantly improved by using a
preconditioner, such as an incomplete factorization, which works well
on general problems, or an application-specific preconditioner, such
as an overlapping additive Schwartz method \cite{petsc}.

In a parallel environment, the cost of message passing magnifies the
importance of an appropriate selection of the linear system solver and
preconditioner.  Parallel implementations of GMRES, conjugate
gradients, and many other iterative methods, along with scalable
preconditioners such as incomplete factorizations and additive
Schwartz methods, are provided in the PETSc toolkit
\cite{petsc,PETSc-user-ref}.  These linear system solvers use the
Message Passing Interface (MPI) \cite{using-mpi} for communication
between processors.

Providing the user with sufficient flexibility in the choice of
linear solver requires careful selection and implementation of the
complementarity algorithms.  In particular, the interface to the
numerical objects, iterative methods and linear algebra, must be
separated from their implementations.  Object-oriented techniques
permit the use of serial or parallel data structures with minimal
changes to the interface and allow new implementations of linear
system solvers and associated linear algebra to be incorporated.  TAO
was designed to enable this flexibility and leverages many of the
existing parallel linear algebra tools, such as HYPRE \cite{hypre},
PETSc, and Trilinos \cite{trilinos-web-page}.  The optimization
methods in TAO, which include the complementarity methods and solvers
to minimize an objective function with bounds on the variables, are 
written to scale on parallel machines and expose many of the algorithmic 
details, such as the linear solver, to the user.  This design philosophy 
assumes that the user has knowledge of an application that can be
used to improve the parallel performance.

\section{Complementarity Algorithms}


The finite-dimensional nonlinear complementarity problem defined by a
given function $F:\Re^n\to\Re^n$ is to calculate an $x^* \in \Re^n$
such that $x^* \geq 0$, $F(x^*) \geq 0$, and $\ip{x^*}{F(x^*)} = 0$.
The algorithms implemented in this study use the current iterate to
define an active set, solve a reduced system to calculate a direction,
and then perform a line search to compute a new iterate that
sufficiently decreases a merit function.

\subsection{Active-Set Semismooth Method}


The semismooth algorithm reformulates a given complementarity problem
as a nonsmooth system of equations satisfying a semismooth property by
using an NCP function.  The resulting nonsmooth system of equations is
solved with Newton's method where an element of the B-subdifferential
plays the role of the Jacobian matrix in the direction calculation.
The active-set semismooth method implemented in TAO is based on this
idea.


Mathematically, a function $\phi:\Re^2\to\Re$ is said to be an
NCP function if $\phi(a,b) = 0$ if and only if $a \geq 0$, $b \geq 0$,
and $ab = 0$.  By defining
\[
   \Phi(x) := \left[\begin{array}{c} 
       \phi(x_1, F_1(x)) \\ 
       \phi(x_2, F_2(x)) \\ 
       \vdots \\
       \phi(x_n, F_n(x)) 
   \end{array}\right]
\]
for any NCP function $\phi$, the nonlinear complementarity problem can
be reformulated as finding an $x^* \in \Re^n$ such that $\Phi(x^*)=0$.
In order to globalize a Newton method for solving this system of equations, the 
merit function
\[
\Psi(x) := \frac{1}{2} \| \Phi(x) \|_2^2
\]
is typically used in a line search.


The Fischer-Burmeister function \cite{fischer:special},
\[
   \phi_{FB}(a,b) := a + b - \sqrt{a^2 + b^2},
\]
is one NCP function, where $\Phi_{FB}$ denotes the reformulation of
the nonlinear complementarity problem implied by using $\phi_{FB}$.
Even though $\Phi_{FB}$ is not continuously differentiable, a
Newton method can still be constructed for finding $\Phi_{FB}(x)=0$.


The function $\Phi_{FB}$ is differentiable almost everywhere and
therefore admits a B-subdifferential \cite{qi:convergence}
\[
   \partial_B \Phi_{FB}(x) := \left\{ H \in \Re^{n\times n} \mid 
      \exists \left\{x^k\right\} \subseteq D_{\Phi_{FB}} \mbox{ with }
      \lim_{x^k\to x} \nabla \Phi_{FB}(x^k) = H \right\},
\]
where $D_{\Phi_{FB}}$ denotes the set of points where $\Phi_{FB}$
is differentiable.  In particular,
\[
   \partial_B \Phi_{FB}(x) \subseteq \left\{ D_a(x) + D_b(x) \nabla F(x) \right\}
\]
for nonnegative diagonal matrices $D_a(x)$ and $D_b(x)$ defined componentwise
as follows \cite{facchinei.soares:merit}:
\begin{enumerate}
\item[(a)] If $\|(x_i, F_i(x))\|_2 > 0$, then
\begin{eqnarray*}
\left[D_a(x)\right]_{i,i} := 1 - \frac{x_i}{\|(x_i, F_i(x))\|_2} \\
\left[D_b(x)\right]_{i,i} := 1 - \frac{F_i(x)}{\|(x_i, F_i(x))\|_2}.
\end{eqnarray*}
\item[(b)] Otherwise
\[
\left(\left[D_a(x)\right]_{i,i},\left[D_b(x)\right]_{i,i}\right) \in
\left\{ (1-\alpha,1-\beta) \mid \|(\alpha,\beta)\|_2 \leq 1 \right\}.
\]
\end{enumerate}
Furthermore, the merit function,
$\Psi_{FB}(x):=\frac{1}{2}\|\Phi_{FB}(x)\|_2^2$ is continuously
differentiable with $\nabla \Psi_{FB}(x) = H^T\Phi_{FB}(x)$ for
any $H \in \partial_B \Phi_{FB}(x)$.


The main computational task in a semismooth Newton method
\cite{deluca.facchinei.ea:semismooth} for solving the nonsmooth system
of equations $\Phi_{FB}(x) = 0$ is to calculate a direction by solving
the linear system of equations
\[
    H^k d^k = -\Phi_{FB}(x^k),
\]
where $H^k$ is any element of $\partial_B \Phi_{FB} (x^k)$.  An Armijo
line search \cite{armijo:minimization} along this direction is then
applied to obtain a new iterate that sufficiently decreases the merit
function.  When the full system is solved to find the Newton
direction, one must use a method suitable for nonsymmetric matrices
because of the row scaling implied by the characterization of the
B-subdifferential.  Alternatively, one can solve a reduced system,
where only a nonnegative diagonal perturbation is made 
to $\nabla F(x)$.  This reduced system will be symmetric 
whenever $\nabla F(x)$ is symmetric, removing the restrictions placed 
on the linear solver for the full-space system.



The reduced system is obtained by selecting active and inactive sets of
constraints.  For a fixed $0 \leq \epsilon < 1$, define
\begin{eqnarray*}
\mathcal{A}(x) := \left\{ i \in \left\{1,\ldots,n\right\} \mid [D_b(x)]_{i,i} \leq \epsilon\right\} \\
\mathcal{I}(x) := \left\{ i \in \left\{1,\ldots,n\right\} \mid [D_b(x)]_{i,i} > \epsilon\right\},
\end{eqnarray*}
where $\mathcal{A}(x)$ denotes the active constraints at $x$ and
$\mathcal{I}(x)$ the inactive constraints.  The characterization of
the B-subdifferential can be used to show that
\[
i \in \mathcal{A}(x) \;\Rightarrow\; x_i \leq \kappa F_i(x)
\]
for some $\kappa > 0$.  The latter characterization is used by the
active-set method in in \cite{facchinei.soares:merit}, where they show
that for all $x$ in a sufficiently small neighborhood of $x^*$, the
strongly active components are correctly identified.

We are now ready to state our active-set semismooth algorithm.
\begin{algorithm}{Active-Set Semismooth Method}
\begin{enumerate}
\item Let $F:\Re^n\to\Re^n$, $x^0 \in \Re^n$, $\epsilon \in [0,1)$, 
$\rho > 0$, $p > 2$, $\beta < 1$, and $\sigma \in (0,\frac{1}{2})$ 
be given.  Set $k = 0$.
\item If $\| \Phi(x^k) \|_2 \leq \mbox{tol}$, then stop.
\item Otherwise, choose $D_a(x^k)$ and $D_b(x^k)$ so that 
\[
D_a(x^k) + D_b(x^k)\nabla F(x^k) \in \partial_B \Phi_{FB}(x^k),
\]
and determine $\mathcal{A}^k := \mathcal{A}(x^k)$ and 
$\mathcal{I}^k := \mathcal{I}(x^k)$.
\item Let
\[
  d^k_{\mathcal{A}^k} = -[D_a(x^k)]_{\mathcal{A}^k,\mathcal{A}^k}^{-1} \Phi_{FB}(x^k)_{\mathcal{A}^k},
\]
and approximately solve the reduced system
\begin{eqnarray*}
\left([D_b(x^k)]_{\mathcal{I}^k,\mathcal{I}^k}^{-1}
      [D_a(x^k)]_{\mathcal{I}^k,\mathcal{I}^k} + 
      [\nabla F(x^k)]_{\mathcal{I}^k,\mathcal{I}^k}\right)
d^k_{\mathcal{I}^k} =  \\
-[D_b(x^k)]_{\mathcal{I}^k,\mathcal{I}^k}^{-1} \Phi_{FB}(x^k)_{\mathcal{I}^k} 
-[\nabla F(x^k)]_{\mathcal{I}^k,\mathcal{A}^k} d^k_{\mathcal{A}^k}
\end{eqnarray*}
to find $d^k_{\mathcal{I}^k}$.
\item If the descent test
\[
\nabla \Psi_{FB}(x^k)^T d^k \leq -\rho\|d^k\|_2^p
\]
is not satisfied, set $d^k = -\nabla \Psi_{FB}(x^k)$.
\item Calculate the smallest $i \in \left\{0, 1, \ldots\right\}$ such that
\[
\Psi_{FB}(x^k + \beta^i d^k) \leq \Psi_{FB}(x^k) + \sigma\beta^k\nabla \Psi_{FB}(x^k)^T d^k,
\]
and set $x^{k+1} = x^k + \beta^i d^k$.
\item Set $k = k + 1$, and go to Step 2.
\end{enumerate}
\end{algorithm}


The advantages of this active-set method are that the linear systems of 
equations solved to calculate the Newton directions are of a reduced size 
and they remain symmetric, positive definite if $\nabla F(x)$ is symmetric,
positive definite.  Therefore, the choice of preconditioner and iterative 
method is not restricted to nonsymmetric methods.

The implementation of the active-set method in TAO uses 
$\rho = 10^{-10}$, $p = 2.1$, $\beta = \frac{1}{2}$, 
and $\sigma = 10^{-4}$.  The value of $\epsilon$ is dynamically 
chosen as
\[
\epsilon(x^k) = \frac{\min\left\{\frac{1}{2} \| \Phi(x^k) \|_2^2, 10^{-2}\right\}}{1 + \| \nabla F(x^k) \|_1}.
\]
This choice forces the active-set identification to go to zero as we
approach a solution to the complementarity problem and deals with
possible scaling problems with the Jacobian of $F$, since small values
for $[D_b]_{i,i}$ can still have an effect on the direction
calculation when the Jacobian is poorly scaled.  For mixed
complementarity problems where we have both finite lower and upper
bounds, the Billups formulation is used \cite{billups:algorithms}.
The same active-set identification technique is used with this
formulation.

\subsection{Reduced-Space Method}

The reduced-space method implemented in TAO also selects an active set
and solves a reduced linear system of equations to calculate a
direction.  The iterates in this method remain within the variable
bounds, while no such guarantee is made with the active-set semismooth
method.  This method is motivated by the simplicity of implementation
in a distributed memory computing environment and the computational
efficiency of similar methods observed on several classes of problems
\cite{dirkse.ferris:crash}.

The active and inactive sets used within the reduced-space method are
defined as
\begin{eqnarray*}
\mathcal{A}(x) := \left\{ i \in \left\{1,\ldots,n\right\} \mid x_i = 0 \mbox{ and } F_i(x) > 0 \right\} \\
\mathcal{I}(x) := \left\{ i \in \left\{1,\ldots,n\right\} \mid x_i > 0 \mbox{ or } F_i(x) \leq 0 \right\}.
\end{eqnarray*}
The active set denotes the variables where the lower bound is active
and the function value can be ignored.  The inactive set contains the
remainder of the variables.  At every iteration of the reduced-space
method a direction is calculated by approximately solving the linear
system equation
\[
     \left[\nabla F(x^k)\right]_{\mathcal{I}^k,\mathcal{I}^k} d_{\mathcal{I}^k}
= -F_{\mathcal{I}^k} (x^k)
\]
and setting $d_{\mathcal{A}^k}$ to zero.

A projected line search is then applied to generate the next 
iterate $x^{k+1}$ such that
\begin{equation} \label{next-y}
x^{k+1} = \pi [ x^k + \alpha d^k ],  
\end{equation}
where $\pi$ is the projection onto the variable bounds.
The step size $\alpha$ is chosen such that
$\| F_{\Omega}(\pi [ x^k + \alpha d^k ] ) \|_2 \leq ( 1 - \sigma \alpha ) \| F_{\Omega}(x) \|_2$,
where $F_{\Omega}(x)$ is defined component-wise
by
\begin{equation} \label{proj-gradient}
 \left[ F_{\Omega}(x) \right] _i = \left\{
\begin{array}{lll}
F_i(x) & \mbox{ if } & x_i > 0 \\
\min \{ F_i(x),0 \} & \mbox{ if } & x_i = 0.
\end{array}
\right.
\end{equation}
The line search tries step lengths $\alpha=\beta^j$ for $\beta \in (0,1)$ and
positive integers $j$ such that $\beta^j > \gamma $.
There are no guarantees the direction calculated is a descent direction, so the
line search can terminate with either a new point of sufficient improvement
or a minimum step length.
The implementation uses $\sigma = 10^{-4}$, $\beta = 0.5$, and
a minimum stepsize of $\gamma = 10^{-12}$ for the parameter choices.

If the line search fails to identify such a 
point, this method performs the same line search in the direction 
$d^k = -F(x^k)$.  The algorithm continues until there is a second failure 
in the line search or a stationary point has been found.

The complete reduced-space algorithm follows.
\begin{algorithm}{Active-Set Reduced Space Method}
\begin{enumerate}
\item Let $F:\Re^n\to\Re^n$ and $x^0 \in \Re_+^n$ be given, and set $k = 0$.
\item If $\| F_{\Omega}(x^k) \|_2 \leq tol $, then stop.
\item Let
\begin{eqnarray*}
\mathcal{A}(x) := \left\{ i \in \left\{1,\ldots,n\right\} \mid x_i = 0 \mbox{ and } F_i(x) > 0 \right\} \\
\mathcal{I}(x) := \left\{ i \in \left\{1,\ldots,n\right\} \mid x_i > 0 \mbox{ or } F_i(x) \leq 0 \right\},
\end{eqnarray*}
set $d_{\mathcal{A}^k} = 0$, and approximately solve the linear system
\[
\left[\nabla F(x^k)\right]_{\mathcal{I}^k,\mathcal{I}^k} d_{\mathcal{I}^k}
= -F_{\mathcal{I}^k} (x^k)
\]
to calculate a direction.
\item If possible, calculate the smallest 
$i \in\left\{0,1,\ldots,\bar{i} \right\}$ such that 
\[
\| F_{\Omega}( \pi [x^k+\beta^i d^k ]) \|_2 \leq
(1 - \sigma \beta^i ) \| F_{\Omega}(x^k) \|_2.
\]
\item Otherwise, set $d^k = -F(x^k)$, and apply the line search.  If no 
such $i$ is found providing sufficient decrease, then stop.
\item Set $x^{k+1} = \pi [x^k+\beta^i d^k]$ and $k = k+1$, and go to step 2.
\end{enumerate}
\end{algorithm}

Although no convergence results have been proven for the reduced-space
method, it has been demonstrated to be effective, especially on
monotone applications.  The reduced matrix retains symmetry and positive 
definiteness when they exist in the Jacobian, permitting the use of
a symmetric linear solver.  Furthermore, this
method lends itself to a parallel implementation because it requires
only a few numerical operations other than a linear solver.

\section{Computational Results}

The complementarity methods were implemented in and are distributed with
the Toolkit for Advanced Optimization.  Two
series of tests were performed on the implementations: one to check
the robustness to verify that the methods work on a significant number
of general problems, and the other to test the parallel performance
and scalability of the methods on some trial infinite-dimensional
variational inequalities.  The latter tests also validate the merits
of our design philosophy by showing that the customization of the linear
system solver does have a significant effect on the overall solution time.
For all of these tests, an upper limit of 100 linear system solves was 
placed on the methods.

\subsection{Robustness}

We ran the implemented methods on the complementarity problems contained
in the MCPLIB collection \cite{dirkse.ferris:mcplib} to demonstrate robustness 
on a diverse set of problems, including many small complementarity 
problems and some nonlinear obstacle problems and optimal control 
problems.  These computational tests were performed on Linux workstation
containing a Pentium 4 processor with a clock speed of 1.8 GHz and 512
MB of RAM.

In order to test the robustness, we ran the methods with LUSOL as the
selected linear system solver.  This LU factorization routine is used
by several optimization and complementarity packages, including MINOS
\cite{murtagh.saunders:minos} and PATH \cite{dirkse.ferris:path}, and
is known to be robust and efficient in a serial environment.  We chose
to use a direct factorization for this test to eliminate possible
problems with the choice of preconditioner, iterative method, and
termination tolerances for the inexact linear system solves.

Tables \ref{T1a} and \ref{T1b} report the number of successes and failures for each
method and model in the test set when using LUSOL.  These results
indicate that the active-set semismooth method, which solves 73.7\% of
the problems, is more robust than the reduced-space method, which
solves 65.5\% of the problems.  The overall conclusion to be drawn is
that both methods solve a significant fraction of these test 
problems.  While these robustness results can be improved by 
introducing heuristics into the implementations, we prefer to use 
the standard methods because they are less complicated and work
well on the targeted infinite-dimensional variational inequalities.


\begin{table}[htbp]
\small
\begin{center}
\caption{Performance of methods on MCPLIB.}
\label{T1a}
\begin{tabular}{|l|cc|cc|}
\hline
\multicolumn{1}{|c|}{}&
\multicolumn{2}{c|}{Active-Set Method} &
\multicolumn{2}{c|}{Reduced-Space Method} \\
\multicolumn{1}{|c|}{Problem}&
\multicolumn{1}{c|}{Successes}&
\multicolumn{1}{c|}{Failures}&
\multicolumn{1}{c|}{Successes}&
\multicolumn{1}{c|}{Failures} \\
\hline
ahn & 1 & 0 & 1 & 0 \\
badfree & 1 & 0 & 1 & 0 \\
baihaung & 1 & 0 & 1 & 0 \\
bert\_oc & 4 & 0 & 4 & 0 \\
bertsekas & 6 & 0 & 3 & 3\\
billups & 0 & 3 & 0 & 3 \\
bishop & 0 & 1 & 0 & 1 \\
bratu & 0 & 1 & 1 & 0 \\
cammcf & 0 & 1 & 0 & 1 \\
cgereg & 2 & 20 & 14 &  8\\
choi & 1 & 0 & 1 & 0 \\
colvdual & 2 & 2 & 2 & 2 \\
colvnlp & 6 & 0 & 6 & 0 \\
cycle & 1 & 0 & 1 & 0 \\
danny & 6 & 2 & 6 & 2 \\
degen & 1 & 0 & 1 & 0 \\
dirkse & 1 & 1 & 1 & 1 \\
duopoly & 0 & 1 & 0 & 1 \\
eckstein & 1 & 0 & 0 & 1 \\
ehl & 2 & 10 & 6 & 6 \\
electric & 0 & 1 & 0 & 1  \\
explcp  & 1 & 0 & 0 & 1  \\
exros & 1 & 4 & 1 & 4 \\
ferrralph & 2 & 0 & 2 & 0 \\
finance & 60 & 0 & 59 & 1 \\
fixedpt & 0 & 2 & 0 & 2 \\
force & 0 & 2 & 0 & 2 \\
freebert & 3 & 4 & 5 & 2 \\
fried & 5 & 5 & 6 & 4 \\
friedms & 1 & 0 & 1 & 0 \\
gafni & 3 & 0 & 3 & 0 \\
games & 18 & 7 & 0 & 25 \\
golanmcp & 0 & 1 & 1 & 0 \\
hanskoop & 7 & 3 & 5 & 5 \\
hydroc & 2 & 0 & 2 & 0 \\
jel & 2 & 0 & 2 & 0 \\
jiangqi & 3 & 0 & 3 & 0 \\
josephy & 8 & 0 & 5 & 3 \\
kanzow & 7 & 0 & 7 & 0 \\
kojshin & 8 & 0 & 5 & 3 \\
kyh & 0 & 4 & 0 & 4 \\
lincont & 0 & 1 & 0 & 1 \\
lstest & 0 & 1 & 0 & 1 \\
leyffer & 1 & 0 & 1 & 0 \\
\hline
\end{tabular}
\end{center}
\end{table}

\begin{table}[htbp]
\small
\begin{center}
\caption{Performance of methods on MCPLIB.}
\label{T1b}
\begin{tabular}{|l|cc|cc|}
\hline
\multicolumn{1}{|c|}{}&
\multicolumn{2}{c|}{Active-Set Method} &
\multicolumn{2}{c|}{Reduced-Space Method} \\
\multicolumn{1}{|c|}{Problem}&
\multicolumn{1}{c|}{Successes}&
\multicolumn{1}{c|}{Failures}&
\multicolumn{1}{c|}{Successes}&
\multicolumn{1}{c|}{Failures} \\
\hline
mathi & 13 & 0 & 13 & 0 \\
methan & 1 & 0 & 1 & 0 \\
mr5mcf & 0 & 1 & 0 & 1 \\
munson & 3 & 1 & 3 & 1 \\
nash & 4 & 0 & 4 & 0 \\
ne-hard & 0 & 1 & 1 & 0 \\
obstacle & 8 & 0 & 8 & 0 \\
optcont & 5 & 0 & 5 & 0 \\
pgvon & 0 & 12 & 2 & 10 \\
pies & 1 & 0 & 0 & 1 \\
pizer & 1 & 3 & 3 & 1 \\
powell & 5 & 1 & 3 & 3 \\
powellmcp & 6 & 0 & 6 & 0 \\
poz & 6 & 0 & 6 & 0 \\
qp & 1 & 0 & 1 & 0 \\
runge & 3 & 4 & 0 & 7 \\
scarf & 9 & 3 & 6 & 6 \\
shansim & 1 & 0 & 1 & 0 \\
shubik & 11 & 37 & 10 & 38 \\
simple-ex & 0 & 1 & 0 & 1 \\
simple-red & 1 & 0 & 1 & 0 \\
spillmcp & 0 & 1 & 0 & 1 \\
sppe & 2 & 1 & 3 & 0 \\
sun & 1 & 0 & 1 & 0 \\
taji & 12 & 0 & 12 & 0 \\
tiebout & 0 & 6 & 0 & 6 \\
tinloi & 64 & 0 & 63 & 1 \\
tinsmall & 63 & 1 & 63 & 1 \\
titan & 1 & 1 & 1 & 1 \\
tobin & 4 & 0 & 4 & 0 \\
tqbilat & 1 & 1 & 2 & 0 \\
trafelas & 0 & 2 & 0 & 2 \\
trig & 2 & 1 & 0 & 3 \\
vonthmcf & 0 & 1 & 0 & 1 \\
xiaohar & 4 & 0 & 3 & 1 \\
xu & 35 & 0 & 5 & 30\\
\hline 
Total & 437 & 156 & 391 & 202 \\
\hline
\end{tabular}
\end{center}
\end{table}

\subsection{Scalability}

To demonstrate the scalability and performance gains that can be made
by selecting appropriate linear system solvers, we tested the
complementarity algorithms implemented on some discretizations of
infinite-dimensional variational inequalities.  All of the
computational tests in this subsection were performed on a Linux
cluster composed of 350 Pentium Xeon processors with a clock speed of
2.4 GHz.  Each node has a minimum of 1 GB of RAM is connected by a
Myrinet 2000 network.

\begin{figure} 
\centerline{\epsfig{figure = 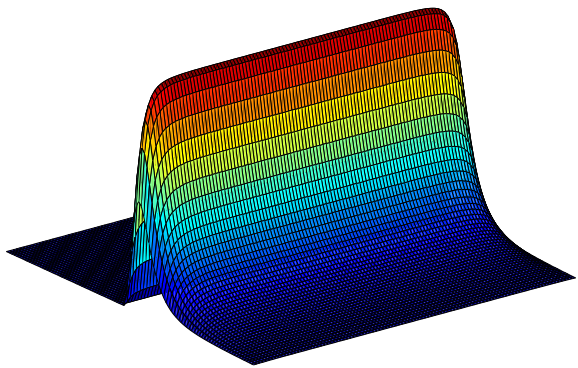, height=2.5in}}
\caption{The journal bearing problem with $ \varepsilon = 0.9 $.\label{pjb}}
\end{figure}

One benchmark application is the journal bearing model, a variational
problem over a two-dimensional region. This problem arises in the
determination of the pressure distribution in a thin film of lubricant
between two circular cylinders.  The infinite-dimensional version of
this problem is to find a piecewise continuously differentiable
function, $v : \cD \mapsto \R $, such that
\[
0 \leq v,\; w_q \Delta v \geq w_l,\; \mbox{ and } v \left[w_q\Delta v - w_l\right] = 0
\]
almost everywhere on $\cD$ with $v = 0$ on $\partial \cD$, where 
$ \cD = ( 0 , 2 \pi ) \times ( 0 , 2b ) $ for some constant 
$ b > 0 $ and
\[
 w_q ( \xi_1 , \xi_2 ) = ( 1 + \varepsilon \cos \xi_1 ) ^ 3 , \quad
 w_l ( \xi_1 , \xi_2 ) = \varepsilon \sin \xi_1
\]
with $ \varepsilon $ in $ (0,1) $. The eccentricity parameter,
$\varepsilon$, influences, in particular, the difficulty of the
problem.  Figure \ref{pjb} shows the solution of the journal bearing
problem for $ \varepsilon = 0.9 $. The steep gradient in the solution
makes this problem a difficult benchmark.  Other elliptic problems
tested include the obstacle, elastic-plastic torsion, and combustion
models from MINPACK-2; see \cite{averick.carter.ea:minpack-2} for a
complete description of these models.

Discretization of the journal bearing problem with either finite
differences or finite elements leads to a standard complementarity
problem.
The number of variables is $ n = n_x n_y $, where $ n_x $ and $ n_y $
are, respectively, the number of grid points in each coordinate
direction of the domain $ \cD $.  See for \cite{more.toraldo:solution}
a description of the finite-element discretization.

The number of grid points can be very large.  In
fact, these problems become so large that the Jacobian matrix cannot
even be stored on one machine.  Direct factorizations exacerbate the
computer memory issues because the fill-in associated with a direct
factorization is significant.  For this reason, iterative linear
solvers are necessary.

The Jacobian of this problem has a sparse, symmetric, and positive
definite structure.  Considerable improvement can be made by using the
conjugate gradient method to solve the linear systems with an
incomplete factorization used as a preconditioner.  Many incomplete
factorization preconditioners require nonzeros along the diagonal, and
the journal bearing problem satisfies this requirement.


The problem was discretized into 40,000 variables and then solved with
three different solvers: LUSOL, a conjugate gradient method with an 
incomplete LU preconditioner, and a conjugate gradient method with
a block Jacobi preconditioner such that each block contains an 
incomplete LU factorization.
As Table \ref{T2} shows, use of the conjugate gradient methods with
with an appropriate preconditioner saves a significant amount of time.  With
the reduced-space method, the
savings range from $30\%$ on the combustion problem to over $70\%$ on
the elastic-plastic torsion problem.

\begin{table}[htbp]
\begin{center}
\caption{Solution times for three linear solvers (sec).}
\label{T2}
\begin{tabular}{|l|ccc|}
\multicolumn{1}{l}{Problem}&
\multicolumn{1}{c}{LUSOL}&
\multicolumn{1}{c}{CG - 1 P}&
\multicolumn{1}{c}{CG - 2 P} \\
\hline
J Bearing & 278 & 136 & 99 \\
EP Torsion & 172 & 50 & 34.6 \\
Obstacle & 49 & 18.4 & 14.3 \\
Combustion & 29.8 & 20.2 & 14.1 \\
\hline
\end{tabular}
\end{center}
\end{table}

To demonstrate the parallel efficiency of the complementarity methods,
we wrote parallel implementations of the elastic-plastic torsion and
obstacle problems using the grid management facilities of PETSc, which
relies on MPI \cite{using-mpi} for communications between processors.
PETSc provides support for discretizing the rectangular region
${\cD}$, partitioning the surface into multiple regions, and assigning
each processor to one of these regions.  Each processor computes the
function on its subdomain with respect to the variables in its region.  A
preconditioned conjugate gradient method was used to solve the linear
systems generated by the complementarity algorithms.  The
preconditioner was a block Jacobi preconditioner with incomplete LU
factorization.  In the reduced-space method, the relative tolerance
used during the linear solve was $0.01$.  The results are summarized
in Table \ref{flops} when using 1--64 processors.

\begin{table}[htbp]
\caption{Scalability.}
\label{flops}
\small
\begin{center}
\begin{tabular}{|llcc|ccccccc|}
\hline
\multicolumn{1}{|c}{Problem}&
\multicolumn{1}{c}{Solver}&
\multicolumn{1}{c}{mx}&
\multicolumn{1}{c}{my}&
\multicolumn{7}{c|}{Processors} \\
\hline
  &   &   &      & 1 & 2 & 4 & 8 & 16 & 32  & 64\\
\hline
Obstacle & R-S & 800 & 800 & 2848 & 1722 & 887 & 484 & 252 & 129 & 71 \\
Obstacle & A-S & 800 & 800 & 975  & 582  & 300 & 168 & 87  & 51  & 30 \\
JBearing & R-S & 800 & 800 & 13597 & 7786& 7021 & 3555 & 1861 & 930 & 544 \\
JBearing & A-S & 800 & 800 & 14025 & 7770& 7019 & 3540 & 1973 & 1082 & 648 \\
Eptorsion   & R-S & 800 & 800 & 835   & 548 & 287 & 152 & 112 & 58  & 41 \\
Eptorsion   & A-S & 800 & 800 & 4744  & 3118& 1644& 877 & 463 & 261 & 165 \\
Combustion  & R-S  & 800 & 800 & 332 & 228  & 122 & 63  & 31  & 15 & 8.6  \\
Combustion  & A-S  & 800 & 800 & 339 & 233  & 124 & 67  & 32  & 17 & 10.2   \\
\hline
\end{tabular}
\end{center}
\end{table}

\begin{figure}[ht]
\begin{center}
   \subfigure[Overall Implementation Efficiency]{
   \label{G14}
   \includegraphics[width=0.45\textwidth,height=0.4\textwidth]{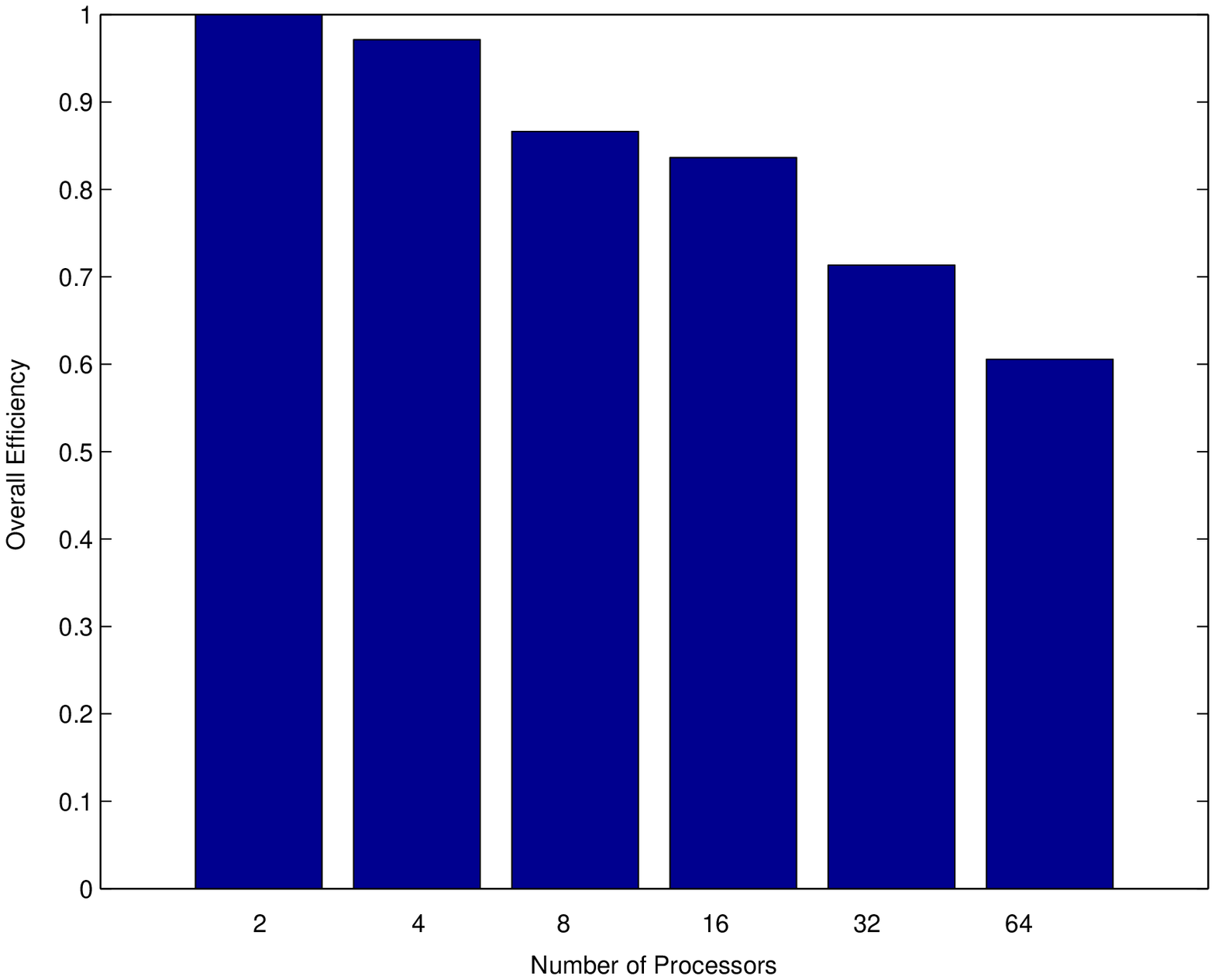}}
   \hskip 0.00\textwidth
   \subfigure[Floating-Point Efficiency]{
   \label{G14chol}
   \includegraphics[width=0.45\textwidth,height=0.4\textwidth]{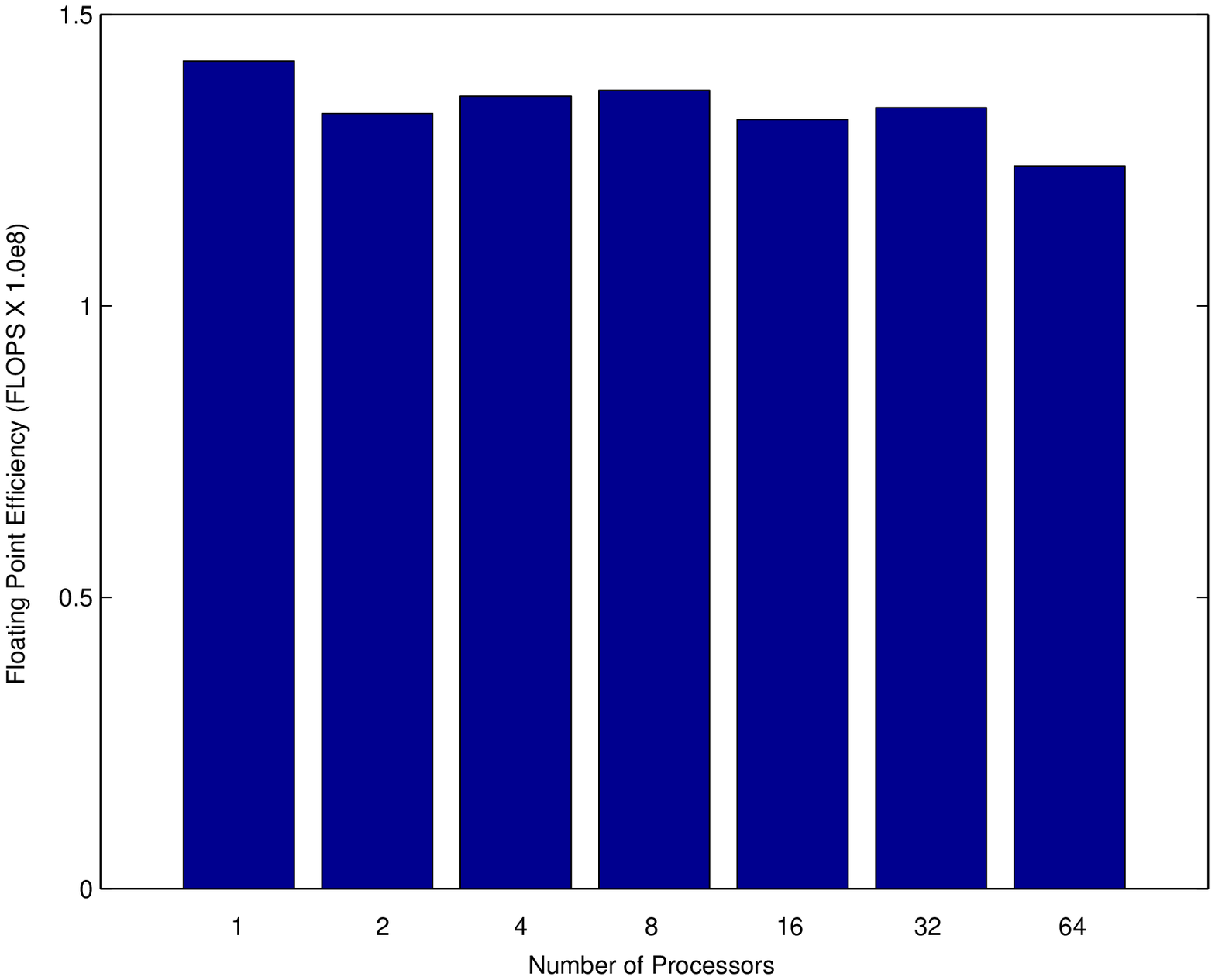}} \\
\caption{Parallel efficiency on the obstacle problem.}
\label{figure2}
\end{center}
\end{figure}

The overall efficiency of our implementation is shown in Figure
\ref{G14}.  Each bar indicates the performance of the active-set 
semismooth method relative
to its performance on two processors.  This number is the ratio of two
times the execution time using two processors and the number of
processors multiplied by the time needed to solve the problem using
one of those processors.  With this metric, the overall parallel efficiency
active-set semismooth method on the obstacle problems using 16 processors was
over $83\%$, and the overall parallel efficiency using 64 processors
was over $60\%$.  The overall efficiency of the reduced-space method 
was similar.

A second set of tests was performed with the mesh refined
according the number of processors used to solve the problem.  In
these tests each processor owned $10,000$ variables.  The mesh was
$100 \times 100$ when one processor was used and $800 \times 800$ when
$64$ processors were used.  Since the methods require more iterations
to solve problems with finer meshes, a comparison of running times
would not demonstrate the scalability of the methods.  Instead, we
used the rate of floating-point operations as the measure of
efficiency.  Figure \ref{G14chol} shows the average number of floating-point 
operations performed per second on each processor when the reduced space method
was used to solve the torsion problem.  As the
overhead of message passing increases, the rate of computation on each
processor decreases.  
The problem used over $142$ MFlops when
one processor was used and that rate only decreased to about  $124$ MFlops when
when 64 processors were used.
Our parallel efficiency of both methods by this measure often
exceeded $80\%$  on sixty four processors.

\section{Conclusion}

The complementarity methods implemented within the Toolkit for
Advanced Optimization are reasonably robust and can effectively
solve discretized versions of infinite-dimensional variational
inequalities.  In order to achieve high parallel performance,
the methods were implemented so that the user has flexibility
in the choice of data structures, linear algebra, and 
algorithms.


\bibliographystyle{plain}
\bibliography{mathprog,tao}

\vskip1in
\begin{center}
\scriptsize
\framebox{\parbox{2.4in}{The submitted manuscript has been created
by the University of Chicago as Operator of Argonne
National Laboratory ("Argonne") under Contract No.\
W-31-109-ENG-38 with the U.S. Department of Energy.
The U.S. Government retains for itself, and others
acting on its behalf, a paid-up, nonexclusive, irrevocable
worldwide license in said article to reproduce,
prepare derivative works, distribute copies to the
public, and perform publicly and display publicly, by or on
behalf of the Government.}}
\normalsize
\end{center}

\end{document}